\documentclass{elsart}
\usepackage{amssymb}
\usepackage{amsfonts}
\usepackage{amsmath}

\input{tcilatex}
\begin{document}

\begin{frontmatter}%

\title{The Hessenberg matrices and Catalan and its generalized numbers}%

\author{Jishe Feng}%

\address{Department of Mathematics, Longdong University,  Qingyang,  Gansu,  745000,  China
E-mail: gsfjs6567@126.com.}%

\begin{abstract}%

We present determinantal representations of the Catalan numbers, $k$%
-Fuss-Catalan numbers, and its generalized number.\ The\ entries of the
normalized Hessenberg matrices are the binomial coefficients that related
with the enumeration of lattice paths.

2010 Mathematics Subject Classification: 05A15, 11B65.

\end{abstract}%

\begin{keyword}%
Hessenberg matrix, Catalan number, $k$-Fuss-Catalan numbers, determinant,
enumeration of lattice paths%
\end{keyword}%

\end{frontmatter}%

\section{Introduction}

An upper Hessenberg matrix $H_{n}=(h_{ij}),i,j=1,2,\cdots ,n$, is a special
kind of square matrix, such that $h_{i,j}=0$ for $i>j+1$. Ulrich Tamm \cite%
{Tamm} give the concept of the Hessenberg matrix in a normalized form, i.e. $%
h_{i+1,i}=1$ for $i=1,\cdots ,n$.\ 
\begin{equation*}
H_{n}=\left( 
\begin{array}{cccccc}
h_{1,1} & h_{1,2} & h_{1,3} & \cdots & \cdots & h_{1,n} \\ 
1 & h_{2,2} & h_{2,3} & \cdots & \cdots & h_{2,n} \\ 
& 1 & h_{3,3} & \cdots & \cdots & h_{3,n} \\ 
&  & \ddots & \ddots &  &  \\ 
&  &  & \ddots & \ddots & h_{n-1,n} \\ 
&  &  &  & 1 & h_{n,n}%
\end{array}%
\right) ,
\end{equation*}

The so-called Pascal matrix $P=\left( \binom{i}{j}\right) _{i,j\geq 0}$ is a
triangular array of the binomial coefficients. We can generalize Pascal
matrix to Hessenberg matrix which elements are binomial coefficients. It is
known that there are a lot of relations between determinants of matrices and
well-known number sequences (see\cite{Feng} \cite{Yilmaz} and references
therein). In this paper, we give a normalized Hessenberg matrices
representation of the famous Catalan and its generalized number.

Lattice paths are omnipresent in enumerative combinatorics, since they can
represent a plethora of different objects. Especially, lattice paths from $%
(0,0)$ to $(x,y)$ with $E=(1,0)$ step and $N=(0,1)$ step that never go above
the line $L:$ $y=kx$, have been models in many combinatorial problems. Let $%
n $ be a positive integer. It is well-known that when $k$ is a positive
integer, the number of lattice paths from $(0,0)$ to $(n,kn)$ which may
touch but never rise above $L$ is $\frac{1}{kn+1}\binom{(k+1)n}{n}$. In
particular, when $k=1$, the number of lattice paths of length $n$ is the $%
\mathit{n}$th Catalan number%
\begin{equation*}
C_{n}=\frac{1}{n+1}\binom{2n}{n}\text{, for n}\geq 0.
\end{equation*}%
The first Catalan numbers for $n=0,1,2,3,...$ are

1, 1, 2, 5, 14, 42, 132, 429, 1430, 4862, 16796, ... (sequence A000108 in
the OEIS).

In general for positive integer $k$, in paper \cite{An}, it is called as $k$%
-Fuss-Catalan numbers $F_{n}$. The purpose of this paper is that we will
give a determinant of normalized Hessenberg matrix in Section 2 for positive
integer, and in Section 3 for rational $k$. In Section 4, we give the
iterative method to evaluate the determinant. Subsequent reduction of
Hessenberg matrix to a triangular matrix can be achieved through iterative
procedures, this is a fast method to evaluate the determinant of upper
normalized Hessenberg matrix.

In \cite{Krattenthaler}, it is to count paths in a region that is delimited
by nonlinear upper and lower boundaries. Let $a_{1}\leq a_{2}\leq \cdots
\leq a_{n}$, and $b_{1}\leq b_{2}\leq \cdots \leq b_{n}$ be integers with $%
a_{i}\geq b_{i}$. We abbreviate $\mathbf{a}=(a_{1},a_{2},...,a_{n})$ and $%
\mathbf{b}=(b_{1},b_{2},...,b_{n})$. Let $L(0,b_{1})\rightarrow (n,a_{n})$
denote the set of all lattice paths from $(0,b_{1})$ to $(n,a_{n})$
satisfying the property that for all $i=1,2,...,n$ the height of the $i$-th
horizontal step is in the interval $[b_{i},a_{i}]$. Theorem 10.7.1 in \cite%
{Krattenthaler} gives a formula for counting these paths, we restate it as
follows.

\begin{thm}
\cite{Krattenthaler} Let $\mathbf{a}=(a_{1},a_{2},...,a_{n})$ and $\mathbf{b}%
=(b_{1},b_{2},...,b_{n})$ be integer sequences with $a_{1}\leq a_{2}\leq
\cdots \leq a_{n}$, $b_{1}\leq b_{2}\leq \cdots \leq b_{n}$, and $a_{i}\geq
b_{i}$, $i=1,2,...,n$. The number of all paths from $(0,b_{1})$ to $%
(n,a_{n}) $ satisfying the property that for all $i=1,2,...,n$ the height of
the $i$-th horizontal step is between $b_{i}$ and $a_{i}$ is given by%
\begin{equation*}
|L((0,b_{1})\rightarrow (n,a_{n}):\mathbf{b}\leq y\leq \mathbf{a})|=\underset%
{1\leq i,j\leq n}{det}\left( \binom{a_{i}-b_{j}+1}{j-i+1}\right) .
\end{equation*}
\end{thm}

\section{The $k$-Fuss-Catalan numbers}

Given a positive number $k$, a $k$-Fuss-Catalan path of length $n$ is a path
from ${(0,0)}$ to $(n,kn)$ using east steps $(1,0)$ and north steps $(0,1)$
such that it stays weakly below the line $y=kx$. The number of all $k$%
-Fuss-Catalan paths of length $n$ is given by the $k$-Fuss-Catalan numbers 
\cite{An},%
\begin{equation*}
F_{n}=\frac{1}{kn+1}\binom{(k+1)n}{n},
\end{equation*}%
and Armstrong\cite{Drew} enumerates the number of $k$-Fuss-Catalan paths of
given type. We can give another formula to $F_{n}$ by the following theorem.

\begin{thm}
Suppose $k$, $n$ are positive integers. Let%
\begin{equation*}
a_{i}=k(i-1)\text{, \ for }i=1,\cdots ,n.
\end{equation*}%
Then the $k$-Fuss-Catalan numbers $F_{n}$ is equal to%
\begin{equation*}
F_{n}=\frac{1}{kn+1}\binom{(k+1)n}{n}=\underset{1\leq i,j\leq n}{det}\left( 
\binom{a_{i}+1}{j-i+1}\right) .
\end{equation*}
\end{thm}

\begin{pf}
Let $a_{i}=k(i-1)$, \ for $i=1,\cdots ,n$, and let 
\begin{equation*}
\mathbf{a}=(a_{1},a_{2},...,a_{n})=(k,2k,\cdots ,k(n-1))
\end{equation*}%
and $\mathbf{b}=(0,0,\cdots 0)$. Applying Theorem 1, we find the number of
lattice paths is $\underset{1\leq i,j\leq n}{det}\left( \binom{a_{i}+1}{j-i+1%
}\right) $.
\end{pf}

It might be interesting to have a formula for enumerating Fuss-Catalan paths
using Theorem 5. We implement an evaluation by a procedure of Maple as
follows.

\textit{\TEXTsymbol{>}with(combinat); \newline
\ \ \ \ with(LinearAlgebra);}\newline
\ \ \ \textit{ffT := proc (n, r)}\newline
\textit{\ \ \ \ \ \ \ \ \ local a, k, j, i, A;}\newline
\ \ \ \ \ \ \ \ \ \textit{k := r*n; \newline
\ \ \ \ \ \ \ \ \ for j to k do a[j] := r*(j-1) end do;\newline
\ \ \ \ \ \ \ \ A := Matrix(n, n, (i, j) -\TEXTsymbol{>}binomial(a[i]+1,
j-i+1));}\newline
\ \ \ \ \ \textit{\ \ \ [A, Determinant(A)];}\newline
\ \ \textit{end proc;}

For $k=3$, we can illustrate the identity by a Maple command:

\textit{\TEXTsymbol{>}for i from 1 to 18 do} \textit{ffT(i,
3)[2];binomial((3+1)*i, i)/(3*i+1);} \textit{end do;}

$1,4,22,140,969,7084,53820,420732,3362260,27343888,\cdots $

For $k=2$, one can get Ternary number \cite{Oeis},%
\begin{equation*}
T_{n}=\frac{1}{2n+1}\binom{3n}{n}=\left\vert 
\begin{array}{cccccc}
\binom{1}{1} &  &  &  &  &  \\ 
\binom{3}{0} & \binom{3}{1} & \binom{3}{2} & \binom{3}{3} &  &  \\ 
& \binom{5}{0} & \binom{5}{1} & \binom{5}{2} & \binom{5}{3} & \cdots \\ 
&  & \ddots & \ddots & \ddots &  \\ 
&  &  & \ddots & \ddots & \ddots \\ 
&  &  &  & \binom{2n-1}{0} & \binom{2n-1}{1}%
\end{array}%
\right\vert .
\end{equation*}%
To illustrate the identity, we use the following Maple command:

\textit{\TEXTsymbol{>}for i from 1 to 18 do ffT(i, 2)[2];binomial((2+1)*i,
i)/(2*i+1);} \textit{end do;}

$1,3,12,55,273,1428,7752,43263,246675,1430715,\cdots $

In particular, when $k=1$, we get an identity about Catalan number,%
\begin{equation*}
C_{n}=\frac{1}{n+1}\binom{2n}{n}=\left\vert 
\begin{array}{cccccc}
\binom{1}{1} &  &  &  &  &  \\ 
\binom{2}{0} & \binom{2}{1} & \binom{2}{2} &  &  &  \\ 
& \binom{3}{0} & \binom{3}{1} & \binom{3}{2} & \binom{3}{2} &  \\ 
&  & \ddots & \ddots & \ddots &  \\ 
&  &  & \ddots & \ddots & \ddots \\ 
&  &  &  & \binom{n}{0} & \binom{n}{1}%
\end{array}%
\right\vert .
\end{equation*}%
\newline

\textit{\TEXTsymbol{>}for i from 1 to 18 do ffT(i, 1)[2];binomial((2+1)*i,
i)/(2*i+1);} \textit{end do;}

Namely, there are%
\begin{eqnarray*}
C_{1} &=&\left\vert 1\right\vert =1;C_{2}=\left\vert 
\begin{array}{cc}
\binom{1}{1} &  \\ 
\binom{2}{0} & \binom{2}{1}%
\end{array}%
\right\vert =2; \\
C_{3} &=&\left\vert 
\begin{array}{ccc}
\binom{1}{1} &  &  \\ 
\binom{2}{0} & \binom{2}{1} & \binom{2}{2} \\ 
& \binom{3}{0} & \binom{3}{1}%
\end{array}%
\right\vert =5;C_{4}=\left\vert 
\begin{array}{cccc}
\binom{1}{1} &  &  &  \\ 
\binom{2}{0} & \binom{2}{1} & \binom{2}{2} &  \\ 
& \binom{3}{0} & \binom{3}{1} & \binom{3}{2} \\ 
&  & \binom{4}{0} & \binom{4}{1}%
\end{array}%
\right\vert =14; \\
&&\cdots
\end{eqnarray*}

\section{The generalized Fuss-Catalan numbers}

When $k=\frac{r}{m}$ is rational, here $r$ and $m$ are coprime positive
integers, it is shown in \cite{Bizley} that the number of lattice paths from 
$(0,0)$ to $(mn,rn)$ that may touch but never rise above the line $y=\frac{r%
}{m}x$, is $\tsum\limits_{a}$ $\tprod\limits_{a_{i}}\frac{F_{i}^{a_{i}}}{%
\alpha _{i}!}$, where $F_{i}=\frac{1}{i(m+r)}\binom{i(m+r)}{im}$ and the sum 
$\tsum\limits_{a}$ is taken over all sequences of non-negative integers $%
a=(a_{1},a_{2},\cdots )$ such that $\tsum\limits_{i=1}^{\infty }ia_{i}=n$.
We give another representation of determinant of normalized upper Hessenberg
matrix as the follows.

\begin{thm}
Suppose $m$, $r$ are coprime positive integers. Let%
\begin{equation*}
a_{i}=r\left\lfloor \frac{i-1}{m}\right\rfloor +\left\lfloor \frac{r}{m}%
(i-m\left\lfloor \frac{i-1}{m}\right\rfloor -1)\right\rfloor \text{, \ for }%
i=1,\cdots ,mn.
\end{equation*}%
Then the generalized Fuss-Catalan numbers $W_{n}$, which is the number of
lattice paths from $(0,0)$ to $(mn,rn)$ that may touch but never rise above
the line $y=\frac{r}{m}x$, is equal to%
\begin{equation*}
W_{n}=\underset{1\leq i,j\leq mn}{det}\left( \binom{a_{i}+1}{j-i+1}\right) .
\end{equation*}

\begin{pf}
Let $a_{i}=r\left\lfloor \frac{i-1}{m}\right\rfloor +\left\lfloor \frac{r}{m}%
(i-m\left\lfloor \frac{i-1}{m}\right\rfloor -1)\right\rfloor $, \ for $%
i=1,\cdots ,nm$, and let $\mathbf{a}=(a_{1},a_{2},...,a_{n})$ and $\mathbf{b}%
=(0,0,\cdots 0)$. Applying Theorem 1, we find the number of lattice paths is 
$\underset{1\leq i,j\leq n}{det}\left( \binom{a_{i}+1}{j-i+1}\right) $.
\end{pf}
\end{thm}

For example, the generalized Fuss-Catalan numbers $W_{1}$, which is the
number of lattice paths from $(0,0)$ to $(7,16)$ that may touch but never
rise above the line $y=\frac{16}{7}x$, is equal to%
\begin{equation*}
W_{1}=\underset{1\leq i,j\leq 7}{det}\left( 
\begin{array}{ccccccc}
\binom{1}{1} &  &  &  &  &  &  \\ 
\binom{3}{0} & \binom{3}{1} & \binom{3}{2} & \binom{3}{3} &  &  &  \\ 
& \binom{5}{0} & \binom{5}{1} & \binom{5}{2} & \binom{5}{3} & \binom{5}{4} & 
\binom{5}{5} \\ 
&  & \binom{7}{0} & \binom{7}{1} & \binom{7}{2} & \binom{7}{3} & \binom{7}{4}
\\ 
&  &  & \binom{10}{0} & \binom{10}{1} & \binom{10}{2} & \binom{10}{3} \\ 
&  &  &  & \binom{12}{0} & \binom{12}{1} & \binom{12}{2} \\ 
&  &  &  &  & \binom{14}{0} & \binom{14}{1}%
\end{array}%
\right) =10659.
\end{equation*}

For the general, we implement the evaluation of the determinant of
Hessenberg matrices by a procedure of Maple as follows.

\TEXTsymbol{>}\textit{with(combinat); \newline
\ \ \ with(LinearAlgebra);\newline
\ \ FjsDyckpath := proc (m, r, n)\newline
\ \ \ local c, k, j, i, C; \newline
\ \ \ \ \ k := m*n; \newline
\ \ \ \ \ for j to k do \newline
\ \ \ \ \ \ \ if j \TEXTsymbol{<}= m then \newline
\ \ \ \ \ \ \ \ \ \ c[j] := floor(r*(j-1)/m) \newline
\ \ \ \ \ \ \ else \newline
\ \ \ \ \ \ \ \ \ \ c[j] :=
floor(r*(j-m*floor((j-1)/m)-1)/m)+r*floor((j-1)/m);\newline
} \ \ \ \ \ \textit{end if }

\textit{\ \ end do;\newline
\ \ \ C := Matrix(m*n, m*n, (i, j)-\TEXTsymbol{>} binomial(c[i]+1, j-i+1); 
\newline
\ \ \ [C, Determinant(C)];\newline
\ end proc}

\TEXTsymbol{>}\textit{for i to 16 do FjsDyckpath(7, 16, i) end do}

\section{Evaluate the determinants}

We can find that above determinants are of upper normalized Hessenberg
matrices, which entries are binomial coefficients. The row operations may be
applied to the matrix {from the first row to the last} row: adding the
negative reciprocal multiple of the entry $a_{ii}$ of the $i$th row to the
succeeding row. This can reduce the matrix to strictly upper triangular, so
we evaluate the determinant of the corresponding upper normalized Hessenberg
matrices. This can reduce the matrix to a strictly upper triangular, so the
determinant of upper normalized Hessenberg matrices is equal to the product
of main diagonal elements.

For example, 
\begin{equation*}
\det \left( 
\begin{array}{cccc}
\binom{1}{1} &  &  &  \\ 
\binom{1}{0} & \binom{1}{1} &  &  \\ 
& \binom{2}{0} & \binom{2}{1} & \binom{2}{2} \\ 
&  & \binom{2}{0} & \binom{2}{1}%
\end{array}%
\right) =\det \left( 
\begin{array}{cccc}
1 &  &  &  \\ 
& 1 &  &  \\ 
&  & 2 & 1 \\ 
&  &  & \frac{3}{2}%
\end{array}%
\right) =3,
\end{equation*}

\begin{equation*}
\det \left( 
\begin{array}{ccccccccc}
\binom{1}{1} &  &  &  &  &  &  &  &  \\ 
\binom{1}{0} & \binom{1}{1} &  &  &  &  &  &  &  \\ 
& \binom{1}{0} & \binom{1}{1} &  &  &  &  &  &  \\ 
&  & \binom{2}{0} & \binom{2}{1} & \binom{2}{2} &  &  &  &  \\ 
&  &  & \binom{2}{0} & \binom{2}{1} & \binom{2}{2} &  &  &  \\ 
&  &  &  & \binom{3}{0} & \binom{3}{1} & \binom{3}{2} & \binom{3}{3} &  \\ 
&  &  &  &  & \binom{3}{0} & \binom{3}{1} & \binom{3}{2} & \binom{3}{3} \\ 
&  &  &  &  &  & \binom{3}{0} & \binom{3}{1} & \binom{3}{2} \\ 
&  &  &  &  &  &  & \binom{4}{0} & \binom{4}{1}%
\end{array}%
\right) =\det \left( 
\begin{array}{ccccccccc}
1 &  &  &  &  &  &  &  &  \\ 
& 1 &  &  &  &  &  &  &  \\ 
&  & 1 &  &  &  &  &  &  \\ 
&  &  & 2 & 1 &  &  &  &  \\ 
&  &  &  & \frac{3}{2} & 1 &  &  &  \\ 
&  &  &  &  & \frac{7}{3} & 3 & 1 &  \\ 
&  &  &  &  &  & \frac{12}{7} & \frac{18}{7} & 1 \\ 
&  &  &  &  &  &  & \frac{18}{12} & \frac{29}{12} \\ 
&  &  &  &  &  &  &  & \frac{43}{18}%
\end{array}%
\right) =43
\end{equation*}

\end{document}